\documentclass[aop,MSNbibl,dvips]{arximspdf}

\doi{10.1214/12-AOP786}
\volume{41}
\issue{3B}
\pubyear{2013}
\firstpage{2090}
\lastpage{2102}

\makeatletter

\newcommand{\eqref}[1]{(\ref{#1})}
\newtheorem{thmm}{Theorem}
\newproclaim{ex}[thmm]{Example}
\makeatother

\begin{document}
\begin{frontmatter}

\title{Exchangeable sequences driven by an absolutely continuous random measure}
\runtitle{Absolutely continuous random measures}

\begin{aug}
\author[A]{\fnms{Patrizia} \snm{Berti}\ead[label=e1]{patrizia.berti@unimore.it}},
\author[B]{\fnms{Luca} \snm{Pratelli}\ead[label=e2]{pratel@mail.dm.unipi.it}}
\and
\author[C]{\fnms{Pietro} \snm{Rigo}\corref{}\ead[label=e3]{pietro.rigo@unipv.it}}
\runauthor{P. Berti, L. Pratelli and P. Rigo}
\affiliation{Universita' di Modena e Reggio-Emilia, Accademia Navale di Livorno
and Universita'~di~Pavia}
\address[A]{P. Berti\\
Dipartimento di Matematica Pura\\
\quad ed Applicata ``G. Vitali''\\
Universita' di Modena e Reggio-Emilia\\
via Campi 213/B\\
41100 Modena\\
Italy\\
\printead{e1}} 
\address[B]{L. Pratelli\\
Accademia Navale di Livorno\\
viale Italia 72\\
57100 Livorno\\
Italy\\
\printead{e2}}
\address[C]{P. Rigo\\
Dipartimento di Matematica ``F. Casorati''\\
Universita' di Pavia\\
via Ferrata 1\\
27100 Pavia\\
Italy\\
\printead{e3}}
\end{aug}

\received{\smonth{3} \syear{2011}}
\revised{\smonth{4} \syear{2012}}

%
\begin{abstract}
Let $S$ be a Polish space and $(X_n\dvtx n\geq1)$ an exchangeable sequence
of $S$-valued random variables.\vspace*{-1pt} Let $\alpha_n(\cdot)=P
(X_{n+1}\in
\cdot\mid X_1,\ldots,X_n )$ be the predictive measure and
$\alpha$
a random probability measure on $S$ such that $\alpha_n\stackrel
{\mathrm{weak}}\longrightarrow\alpha$ a.s. Two (related) problems are
addressed. One is to give conditions for $\alpha\ll\lambda$ a.s., where
$\lambda$ is a (nonrandom) $\sigma$-finite Borel measure on $S$. Such
conditions should concern the finite dimensional distributions
$\mathcal
{L}(X_1,\ldots,X_n)$, $n\geq1$, only. The other problem is to
investigate whether $\Vert\alpha_n-\alpha\Vert\stackrel
{\mathrm{a.s.}}\longrightarrow0$, where $\Vert\cdot\Vert$ is total variation norm.
Various results are obtained. Some of them do not require
exchangeability, but hold under the weaker assumption that $(X_n)$ is
conditionally identically distributed, in the sense of [\textit{Ann. Probab.} \textbf{32} (2004) 2029--2052].
\end{abstract}

%

\end{frontmatter}

\section{Two related problems}\label{intro}

Throughout, $S$ is a Polish space and
\[
X=(X_1,X_2,\ldots)
\]
a sequence of $S$-valued random variables on the probability space
$(\Omega,\mathcal{A},P)$. We let $\mathcal{B}$ denote the Borel
$\sigma$-field on $S$ and $\mathbb{S}$ the set of probability
measures on
$\mathcal{B}$. A~random probability measure on $S$ is a map
$\alpha\dvtx \Omega\rightarrow\mathbb{S}$ such that $\sigma(\alpha
)\subset
\mathcal{A}$, where $\sigma(\alpha)$ is the $\sigma$-field on
$\Omega$
generated by $\omega\mapsto\alpha(\omega)(B)$ for all $B\in
\mathcal{B}$.

For each $n\geq1$, let $\alpha_n$ be the $n$th \textit{predictive
measure}. Thus, $\alpha_n$ is a random probability measure on $S$, and
$\alpha_n(\cdot)(B)$ is a version of $P (X_{n+1}\in B\mid
X_1,\ldots
,X_n )$ for all $B\in\mathcal{B}$. Define also $\alpha_0(\cdot
)=P(X_1\in\cdot)$.

If $X$ is \textit{exchangeable}, as assumed in this section, there
is a random probability measure $\alpha$ on $S$ such that
\[
\alpha_n(\omega)\stackrel{\mathrm{weak}}\longrightarrow\alpha(\omega) \qquad\mbox
{for almost all }\omega\in\Omega.
\]
Such an $\alpha$ can also be viewed as
\[
\mu_n(\omega)\stackrel{\mathrm{weak}}\longrightarrow\alpha(\omega) \qquad\mbox{for
almost all }\omega\in\Omega,
\]
where $\mu_n=\frac{1}{n}\sum_{i=1}^n\delta_{X_i}$ is the empirical
measure. Further, $\alpha$ grants the usual representation
\begin{eqnarray*}
P(X\in B)=\int\alpha(\omega)^\infty(B) P(d\omega)\qquad \mbox{for every Borel
set }B\subset S^\infty,
\end{eqnarray*}
where $\alpha(\omega)^\infty=\alpha(\omega)\times\alpha(\omega
)\times
\cdots$.

Let $\lambda$ be a $\sigma$-finite measure on $\mathcal{B}$. Our
\textit{first problem} is to give conditions for
%
\begin{equation}
\label{goal} \alpha(\omega)\ll\lambda \qquad\mbox{for almost all }\omega\in \Omega.
\end{equation}
The conditions should concern the finite dimensional distributions
$\mathcal{L}(X_1,\ldots,\break X_n)$, $n\geq1$, only.

While investigating \eqref{goal}, one meets another problem, of
possible independent interest. Let $\Vert\cdot\Vert$ denote total
variation norm on $(S,\mathcal{B})$. Our \textit{second problem} is to give
conditions for
\[
\Vert\alpha_n-\alpha\Vert\stackrel{\mathrm{a.s.}}\longrightarrow0.
\]

\section{Motivations}\label{mot}

Again, let $X=(X_1,X_2,\ldots)$ be exchangeable.

Reasonable conditions for \eqref{goal} look of theoretical interest.
They are of practical interest as well thanks to Bayesian
nonparametrics. In this framework, the starting point is a
prior $\pi$ on $\mathbb{S}$. Since $\pi=P\circ\alpha^{-1}$, condition
\eqref{goal} is equivalent to
\[
\pi \{\nu\in\mathbb{S}\dvtx \nu\ll\lambda \}=1.
\]
This is a basic information for the
subsequent statistical analysis. Roughly speaking, it means that the
``underlying statistical model'' consists of absolutely continuous
laws.

Notwithstanding the significance of \eqref{goal}, however, there is a
growing literature which gets around the first problem of this paper.
Indeed, in a plenty of Bayesian nonparametric problems, condition
\eqref
{goal} is just a crude \textit{assumption} and the prior~$\pi$ is directly
assessed on a set of densities (with respect to $\lambda$). See, for
example,
\cite{GV}~and references therein. Instead, it seems reasonable to get
\eqref{goal} as a consequence of explicit assumptions on the finite
dimensional distributions $\mathcal{L}(X_1,\ldots,X_n)$, \mbox{$n\geq1$}.
From a foundational point of view, in fact, only
assumptions on observable facts make sense. This attitude is strongly
supported by de Finetti, among others. When dealing with the sequence
$X$, the observable facts are events of the type $\{(X_1,\ldots
,X_n)\in
B\}$ for some $n\geq1$ and $B\in\mathcal{B}^n$. This is why, in this
paper, the
conditions for \eqref{goal} are requested to concern
$\mathcal{L}(X_1,\ldots,X_n)$, $n\geq1$, only.

Some references related to the above remarks are \cite{CR} and \cite
{DF84,DF88,DF90,DF04B,DF04A,FLR}. In particular, in~\cite{DF04A} and \cite{DF04B},
Diaconis and Freedman have an exchangeable sequence of\vadjust{\goodbreak} indicators and
give conditions for the mixing measure (i.e., the prior $\pi$) to be
absolutely continuous with respect to Lebesgue measure. The present
paper is much in the spirit of \cite{DF04A} and \cite{DF04B}. The main
difference is that we give conditions for the mixands $\{\alpha(\omega
)\dvtx \omega\in\Omega\}$, and not for the mixing measure~$\pi$, to be
absolutely continuous.

Next, a necessary condition for \eqref{goal} is
%
\begin{equation}
\label{nec} \mathcal{L}(X_1,\ldots,X_n)\ll
\lambda^n\qquad \mbox{for all }n\geq1,
\end{equation}
where $\lambda^n=\lambda\times\cdots\times\lambda$. Condition
\eqref
{nec} clearly involves the finite dimensional distributions only. Thus,
a (natural) question is whether
\eqref{nec} suffices for \eqref{goal} as well.

The answer is yes provided $\alpha$ can be approximated by the
predictive measures $\alpha_n$ in some stronger sense. In fact,
condition \eqref{nec} can be written as
\[
\alpha_n(\omega)\ll\lambda \qquad\mbox{for all }n\geq0\mbox{ and almost
all }\omega\in\Omega.
\]
Hence, if \eqref{nec} holds and $\Vert\alpha_n-\alpha\Vert
\stackrel
{\mathrm{a.s.}}\longrightarrow0$, the set
\[
A=\bigl\{\Vert\alpha_n-\alpha\Vert\rightarrow0\bigr\}\cap\{
\alpha_n\ll \lambda\mbox{ for all }n\geq0\}
\]
has probability 1. And, for each $\omega\in A$, one obtains
\[
\alpha(\omega) (B)=\lim_n\alpha_n(\omega) (B)=0\qquad
\mbox{whenever }B\in\mathcal{B}\mbox{ and }\lambda(B)=0.
\]

Therefore, \eqref{goal} follows from \eqref{nec} and $\Vert\alpha_n-\alpha\Vert\stackrel{\mathrm{a.s.}}\longrightarrow0$. In addition, a martingale
argument implies the converse implication, that is,
\[
\alpha\ll\lambda\mbox{ a.s.}\quad\Longleftrightarrow \quad\Vert\alpha_n-
\alpha\Vert\stackrel {\mathrm{a.s.}}\longrightarrow 0\quad\mbox{and}\quad\mathcal{L}(X_1,
\ldots,X_n)\ll\lambda^n\qquad\mbox{for all }n;
\]
see Theorem \ref{p09nh5t}. Thus, our first problem turns into the
second one.

The question of whether $\Vert\alpha_n-\alpha\Vert\stackrel
{\mathrm{a.s.}}\longrightarrow
0$ is of independent interest. Among other things, it is connected
to Bayesian consistency. Surprisingly, however, this question seems
not answered so far. To the best of our knowledge,
$\Vert\alpha_n-\alpha\Vert\stackrel{\mathrm{a.s.}}\longrightarrow0$ in
every example
known so far. And in fact, for some time, we conjectured that
$\Vert\alpha_n-\alpha\Vert\stackrel{\mathrm{a.s.}}\longrightarrow0$ under condition
\eqref{nec}. But this is not true. As shown in Example~\ref{mainexamp}, when $S=\mathbb{R}$ and $\lambda= $ Lebesgue measure,
it may be that
$\mathcal{L}(X_1,\ldots,X_n)\ll\lambda^n$ for all $n$, and yet
$\alpha$
is singular
continuous a.s. Indeed, the (topological) support of $\alpha(\omega)$
has Hausdorff
dimension 0 for almost all $\omega\in\Omega$.

Thus, \eqref{nec} does not suffice for \eqref{goal}. To get
\eqref{goal}, in addition to \eqref{nec}, one needs some growth
conditions on the conditional densities. We refer to forthcoming
Theorem~\ref{rfg76yu8} for such conditions. Here, we mention a result
on the second problem. Actually, for $\Vert\alpha_n-\alpha\Vert
\stackrel
{\mathrm{a.s.}}\longrightarrow0$, it suffices that
\[
P\bigl\{\omega\dvtx \alpha_c(\omega)\ll\lambda\bigr\}=1,
\]
where $\alpha_c(\omega)$ denotes the continuous part of $\alpha
(\omega
)$; see Theorem \ref{gt67fc4}.

Finally, most results mentioned above do not need
exchangeability of $X$, but the weaker assumption
\[
(X_1,\ldots,X_n,X_{n+2})\sim(X_1,
\ldots,X_n,X_{n+1}) \qquad\mbox{for all }n\geq0.\vadjust{\goodbreak}
\]
Those sequences $X$ satisfying the above condition, investigated
in \cite{BPR04}, are called \textit{conditionally identically
distributed} (c.i.d.).\vspace*{-2pt}

\section{Mixtures of i.i.d. absolutely continuous sequences}

In this section, $\mathcal{G}_0=\{\varnothing,\Omega\}$, $\mathcal
{G}_n=\sigma(X_1,\ldots,X_n)$ for $n\geq1$ and $\mathcal{G}_\infty
=\sigma (\bigcup_n\mathcal{G}_n )$. If $\mu$ is a random
probability measure on $S$, we write $\mu(B)$ to denote the real random
variable $\mu(\cdot)(B)$, $B\in\mathcal{B}$. Similarly, if
$h\dvtx S\rightarrow\mathbb{R}$ is a Borel function, integrable with respect
to $\mu(\omega)$ for almost all $\omega\in\Omega$, we write $\mu
(h)$ to
denote $\int h(x) \mu(\cdot)(dx)$.\vspace*{-2pt}

\subsection{Preliminaries}\label{prel}
Let $X=(X_1,X_2,\ldots)$ be c.i.d., as defined in Section \ref{mot}.
Since $X$ needs not be exchangeable, the representation $P(X\in\cdot
)= \int\alpha(\omega)^\infty(\cdot) P(d\omega)$ can fail for any
$\alpha
$. However, there is a random probability measure $\alpha$ on $S$ such that
%
\begin{equation}
\label{cideq23vg} \sigma(\alpha)\subset\mathcal{G}_\infty \quad\mbox{and}\quad
\alpha_n(B)=E \bigl\{\alpha(B)\mid\mathcal{G}_n \bigr\}\qquad
\mbox{a.s.}
\end{equation}
for all $B\in\mathcal{B}$. In particular, $\alpha_n\stackrel
{\mathrm{weak}}\longrightarrow\alpha$ a.s. Also, letting
\[
\mu_n=\frac{1}{n}\sum_{i=1}^n
\delta_{X_i}
\]
be the empirical measure, one obtains $\mu_n\stackrel
{\mathrm{weak}}\longrightarrow\alpha$ a.s. Such an $\alpha$ is of interest for
one more reason. There is an exchangeable sequence $Y=(Y_1,Y_2,\ldots)$
of $S$-valued random variables on $(\Omega,\mathcal{A},P)$ such that
\[
(X_n,X_{n+1},\ldots)\stackrel{d}\longrightarrow Y \quad\mbox
{and}\quad P (Y\in\cdot )=\int\alpha(\omega)^\infty(\cdot) P(d\omega).
\]
See \cite{BPR04} for details.

We next recall some known facts about vector-valued martingales; see~\cite{N}.
Let $(\mathcal{Z},\Vert\cdot\Vert_*)$ be a separable Banach
space. Also, let $\mathcal{F}=(\mathcal{F}_n)$ be a filtration and
$(Z_n)$ a sequence of $\mathcal{Z}$-valued random variables on
$(\Omega
,\mathcal{A},P)$ such that $E\Vert Z_n\Vert_*<\infty$ for all $n$. Then,
$(Z_n)$ is an $\mathcal{F}$-martingale in case $(\phi(Z_n))$ is an
$\mathcal{F}$-martingale for each linear continuous functional $\phi
\dvtx \mathcal{Z}\rightarrow\mathbb{R}$. If $(Z_n)$ is an $\mathcal
{F}$-martingale, $(\Vert Z_n\Vert_*)$ is a real-valued $\mathcal
{F}$-submartingale. So, Doob's maximal inequality yields
\[
E \Bigl\{\sup_n\Vert Z_n\Vert_*^p \Bigr\}
\leq \biggl(\frac
{p}{p-1} \biggr)^p \sup_nE \bigl\{
\Vert Z_n\Vert_*^p \bigr\} \qquad\mbox{for all }p>1.
\]
The following martingale convergence theorem is available as well. Let
$Z\dvtx\break \Omega\rightarrow\mathcal{Z}$ be $\mathcal{F}_\infty
$-measurable and
such that $E\Vert Z\Vert_*<\infty$, where $\mathcal{F}_\infty
=\sigma
(\bigcup_n\mathcal{F}_n )$. Then, $Z_n\stackrel
{\mathrm{a.s.}}\longrightarrow Z$
provided $\phi(Z_n)=E \{\phi(Z)\mid\mathcal{F}_n \}$ a.s. for
all $n$ and all linear continuous functionals $\phi\dvtx \mathcal
{Z}\rightarrow\mathbb{R}$.\vspace*{-2pt}

\subsection{Results}\label{res}

In the sequel, $\lambda$ is a $\sigma$-finite measure on $\mathcal{B}$.
When $S=\mathbb{R}$, it may be natural to think of $\lambda$ as the
Lebesgue measure, but this is only a particular case.\vadjust{\goodbreak} Indeed, $\lambda$
could be singular continuous or concentrated on any Borel subset. In
addition, $X$ is c.i.d. (in particular, exchangeable), and $\alpha$ is
a random probability measure on $S$ such that $\alpha_n\stackrel
{\mathrm{weak}}\longrightarrow\alpha$ a.s. Equivalently, $\alpha$ can be
obtained as $\mu_n\stackrel{\mathrm{weak}}\longrightarrow\alpha$ a.s. It can (and
will) be assumed $\sigma(\alpha)\subset\mathcal{G}_\infty$.

\begin{thmm}\label{p09nh5t} Suppose $X=(X_1,X_2,\ldots)$ is c.i.d. Then,
$\alpha\ll\lambda$ a.s. if and only if
$ \Vert\alpha_n-\alpha\Vert\stackrel{a.s.}\longrightarrow0$ and
$\mathcal
{L}(X_1,\ldots,X_n)\ll\lambda^n$ for all $n$.
\end{thmm}

\begin{pf} The ``if'' part can be proved exactly as in Section \ref
{mot}. Conversely, suppose $\alpha\ll\lambda$ a.s. It can be assumed
$\alpha(\omega)\ll\lambda$ for all $\omega\in\Omega$. We let
$L_p=L_p(S,\mathcal{B},\lambda)$ for each $1\leq p\leq\infty$.

Let $f\dvtx \Omega\times S\rightarrow[0,\infty)$ be such that $\alpha
(\omega
)(dx)=f(\omega,x) \lambda(dx)$ for all $\omega\in\Omega$. Since
$\mathcal{B}$ is countably generated, $f$ can be taken to be $\mathcal
{A}\otimes\mathcal{B}$-measurable (see~\cite{DM}, V.5.58, page 52)
so that
\[
1=\int1 \,dP=\int\!\!\int f(\omega,x) \lambda(dx) P(d\omega)=\int E \bigl\{ f(\cdot,x)
\bigr\} \lambda(dx).
\]
Thus, given $n\geq0$, $E \{f(\cdot,x)\mid\mathcal{G}_n \}
$ is
well defined for $\lambda$-almost all $x\in S$. Since $X$ is c.i.d.,
condition \eqref{cideq23vg} also implies
\begin{eqnarray*}
\int_B E \bigl\{f(\cdot,x)\mid\mathcal{G}_n
\bigr\} \lambda (dx)&=&E \biggl\{ \int_Bf(\cdot,x)
\lambda(dx)\mid\mathcal{G}_n \biggr\}
\\
&=&E \bigl\{\alpha(B)\mid\mathcal{G}_n \bigr\}=\alpha_n(B)\qquad
\mbox{a.s. for fixed }B\in\mathcal{B}.
\end{eqnarray*}
Since $\mathcal{B}$ is countably generated, the previous equality yields
\begin{eqnarray*}
\alpha_n(\omega) (dx)=E \bigl\{f(\cdot,x)\mid\mathcal{G}_n
\bigr\} (\omega ) \lambda(dx) \qquad\mbox{for almost all }\omega\in\Omega.
\end{eqnarray*}
This proves that $\mathcal{L}(X_1,\ldots,X_n)\ll\lambda^n$ for all $n$.
In particular, up to modifying $\alpha_n$ on a $P$-null set, it can be
assumed $\alpha_n(\omega)(dx)=f_n(\omega,x) \lambda(dx)$ for all
$n\geq0$, all $\omega\in\Omega$, and suitable functions $f_n\dvtx \Omega
\times S\rightarrow[0,\infty)$.

Regard $f, f_n\dvtx \Omega\rightarrow L_1$ as $L_1$-valued random
variables. Then, $f\dvtx \Omega\rightarrow L_1$ is $\mathcal{G}_\infty
$-measurable for $\int h(x) f(\cdot,x) \lambda(dx)=\alpha(h)$ is
$\mathcal{G}_\infty$-measurable for all $h\in L_\infty$. Clearly,
$\Vert f(\omega,\cdot)\Vert_{L_1}=\Vert f_n(\omega,\cdot)\Vert_{L_1}=1$ for all $n$
and $\omega$. Finally, $X$ c.i.d. implies
\begin{eqnarray*}
E \biggl\{\int h(x)f(\cdot,x)\lambda(dx)\mid\mathcal{G}_n \biggr\}
&=&E \bigl\{ \alpha(h)\mid\mathcal{G}_n \bigr\} =\alpha_n(h)
\\
&=&\int h(x)f_n(\cdot,x)\lambda(dx) \qquad\mbox{a.s. for all }h\in
L_\infty.
\end{eqnarray*}
By the martingale convergence theorem (see Section \ref{prel})
$f_n\stackrel{\mathrm{a.s.}}\longrightarrow f$ in the space~$L_1$, that is,
\begin{eqnarray}
\bigl\Vert\alpha_n(\omega)-\alpha(\omega)\bigr\Vert=\frac{1}{2} \int
\bigl\vert f_n(\omega,x)-f(\omega,x)\bigr\vert \lambda(dx)\longrightarrow 0
\nonumber\\
\eqntext{\mbox{for almost all }\omega\in\Omega.\qquad}
\end{eqnarray}
\upqed\end{pf}\eject

In the exchangeable case, the argument of the previous proof yields a
little bit more. Indeed, if $X$ is exchangeable and $\alpha\ll\lambda$
a.s., then
\begin{eqnarray*}
\sup_{B\in\mathcal{B}^k}\bigl\vert P \bigl\{(X_{n+1},\ldots ,X_{n+k})\in
B\mid \mathcal{G}_n \bigr\}-\alpha^k(B)\bigr\vert\stackrel
{\mathrm{a.s.}}\longrightarrow0,
\end{eqnarray*}
where $k\geq1$ is any integer and $\alpha^k=\alpha\times\cdots
\times
\alpha$.

The next result deals with the second problem of Section \ref{intro}.
For each $\nu\in\mathbb{S}$, let $\nu_c$ and $\nu_d$ denote the
continuous and discrete parts of $\nu$, that is, $\nu_d(B)=\sum_{x\in
B}\nu\{x\}$ for all $B\in\mathcal{B}$ and $\nu_c=\nu-\nu_d$.

\begin{thmm}\label{gt67fc4} Suppose $X=(X_1,X_2,\ldots)$ is c.i.d. and
$P\{\omega\dvtx \alpha_c(\omega)\ll\lambda\}=1$. Then, $\Vert\alpha_n-\alpha \Vert\stackrel{{a.s.}}\longrightarrow0$ if and only if
%
\begin{equation}
\label{mo87xqml9}
\begin{tabular}{p{220pt}@{}}
\mbox{there is a set }$A_0\in\mathcal{A}$\mbox{ such
that }$P(A_0)=1$\mbox{ and }
$\alpha_n(\omega)\{x\}\longrightarrow\alpha(\omega)\{x\}$ \mbox{ for
all }$x\in S$\mbox{ and }$\omega\in A_0.$
\end{tabular}
\end{equation}
(Recall that $\mathcal{A}$ denotes the basic $\sigma$-field on
$\Omega
$). Moreover, condition \eqref{mo87xqml9} is automatically true if $X$
is exchangeable, so that $\Vert\alpha_n-\alpha\Vert\stackrel
{{a.s.}}\longrightarrow0$ provided $X$ is exchangeable and $\alpha_c\ll
\lambda$ a.s.
\end{thmm}

\begin{pf}
The ``only if'' part is trivial. Suppose condition \eqref{mo87xqml9}
holds. For each $n\geq0$, take functions $\beta_n$ and $\gamma_n$ on
$\Omega$ such that $\beta_n(\omega)$ and $\gamma_n(\omega)$ are
measures on $\mathcal{B}$ for all $\omega\in\Omega$ and
\[
\beta_n(B)=E \bigl\{\alpha_c(B)\mid\mathcal{G}_n
\bigr\},\qquad \gamma_n(B)=E \bigl\{\alpha_d(B)\mid
\mathcal{G}_n \bigr\}\qquad \mbox{a.s.}
\]
for all $B\in\mathcal{B}$. Since $X$ is c.i.d., condition \eqref
{cideq23vg} yields $\alpha_n=\beta_n+\gamma_n$ a.s.

We first prove $\Vert\beta_n-\alpha_c\Vert\stackrel
{\mathrm{a.s.}}\longrightarrow
0$. It can be assumed $\alpha_c(\omega)\ll\lambda$ for all $\omega
\in
\Omega$, so that $\alpha_c(\omega)(dx)=f(\omega,x) \lambda(dx)$ for
all $\omega\in\Omega$ and some function $f\dvtx \Omega\times
S\rightarrow
[0,\infty)$. For fixed $B\in\mathcal{B}$, arguing as in the proof of
Theorem \ref{p09nh5t}, one has
\[
\beta_n(B)=E \biggl\{\int_Bf(\cdot,x)
\lambda(dx)\mid\mathcal {G}_n \biggr\} =\int_BE
\bigl(f(\cdot,x)\mid\mathcal{G}_n \bigr) \lambda (dx)\qquad \mbox{a.s.}
\]
By standard arguments, it follows that $\beta_n\ll\lambda$ a.s. Again,
it can be assumed $\beta_n(\omega)(dx)=f_n(\omega,x) \lambda(dx)$ for
all $\omega\in\Omega$ and some function $f_n\dvtx \Omega\times
S\rightarrow
[0,\infty)$. Define $L_1=L_1(S,\mathcal{B},\lambda)$ and regard
$f_n,
f\dvtx \Omega\rightarrow L_1$ as $L_1$-valued random variables. By the same
martingale argument used for Theorem \ref{p09nh5t}, one obtains
$f_n\stackrel{\mathrm{a.s.}}\longrightarrow f$ in the space $L_1$. That is,
$\Vert\beta_n-\alpha_c\Vert\stackrel{\mathrm{a.s.}}\longrightarrow0$.

We next prove $\Vert\gamma_n-\alpha_d\Vert\stackrel
{\mathrm{a.s.}}\longrightarrow
0$. Take $A_0$ as in condition \eqref{mo87xqml9}, and define
\[
A_1=\Bigl\{ \lim_n \Vert f_n-f
\Vert_{L_1}=0 \mbox{ and } \alpha_n=\beta_n+
\gamma_n\mbox{ for all }n\geq0\Bigr\}.
\]
Then, $P(A_0\cap A_1)=1$ and
\begin{eqnarray*}
\alpha_d(\omega)\{x\}&=&\alpha(\omega)\{x\}-\alpha_c(
\omega)\{x\} =\alpha (\omega)\{x\}-f(\omega,x) \lambda\{x\}
\\
&=&\lim_n \bigl(\alpha_n(\omega)\{x\}-f_n(
\omega,x) \lambda\{x\} \bigr)=\lim_n \bigl(\alpha_n(
\omega)\{x\}-\beta_n(\omega)\{x\} \bigr)\\
&=&\lim_n
\gamma_n(\omega)\{x\}
\end{eqnarray*}
for all $\omega\in A_0\cap A_1$ and $x\in S$. Define also
\[
A=A_0\cap A_1\cap\bigl\{\gamma_n(S)
\longrightarrow\alpha_d(S)\bigr\}.
\]
Since $\gamma_n(S)=1-\beta_n(S)\stackrel{\mathrm{a.s.}}\longrightarrow
1-\alpha_c(S)=\alpha_d(S)$, then $P(A)=1$. Fix $\omega\in A$ and let
$D_\omega
=\{x\in S\dvtx \alpha(\omega)\{x\}>0\}$. Then
\[
\alpha_d(\omega) (D_\omega)\leq\liminf_n
\gamma_n(\omega) (D_\omega)
\]
since $D_\omega$ is countable and $\alpha_d(\omega)\{x\}=\lim_n\gamma_n(\omega)\{x\}$ for all $x\in D_\omega$. Further,
\begin{eqnarray*}
\limsup_n\gamma_n(\omega) (D_\omega)\leq
\limsup_n\gamma_n(\omega ) (S)=\alpha_d(
\omega) (S)=\alpha_d(\omega) (D_\omega).
\end{eqnarray*}
Therefore, $\lim_n\Vert\gamma_n(\omega)-\alpha_d(\omega)\Vert=0$
is an
immediate consequence of
\begin{eqnarray*}
\gamma_n(\omega)\{x\}&\longrightarrow&\alpha_d(\omega)\{x
\} \qquad\mbox{for each }x\in D_\omega,
\\
\alpha_d(\omega) (D_\omega)\hspace*{-3pt}&=&\hspace*{-3pt}\lim_n
\gamma_n(\omega) (D_\omega ), \qquad \alpha_d(\omega)
\bigl(D_\omega^c\bigr)=\lim_n\gamma_n(
\omega) \bigl(D_\omega^c\bigr)=0.
\end{eqnarray*}

Finally, suppose $X$ is exchangeable. We have to prove condition \eqref
{mo87xqml9}. If $S$ is countable, condition \eqref{mo87xqml9} is
trivial for $\alpha_n(B)\stackrel{\mathrm{a.s.}}\longrightarrow\alpha(B)$ for
fixed $B\in\mathcal{B}$. If $S=\mathbb{R}$, the Glivenko--Cantelli
theorem yields $\sup_x\vert\mu_n(I_x)-\alpha(I_x)\vert\stackrel
{\mathrm{a.s.}}\longrightarrow0$, where $I_x=(-\infty,x]$ and $\mu_n=\frac
{1}{n}\sum_{i=1}^n\delta_{X_i}$ is the empirical measure. Hence,
\eqref
{mo87xqml9} follows from
\begin{eqnarray*}
\sup_x \bigl\vert\alpha_n(I_x)-
\mu_n(I_x)\bigr\vert\stackrel {\mathrm{a.s.}}\longrightarrow0;
\end{eqnarray*}
see Corollary 3.2 of \cite{BMR}. If $S$ is any uncountable Polish
space, take a Borel isomorphism $\psi\dvtx S\rightarrow\mathbb{R}$. (Thus
$\psi$ is bijective with $\psi$ and $\psi^{-1}$ Borel measurable). Then
$(\psi(X_n))$ is an exchangeable sequence of real random variables, and
condition \eqref{mo87xqml9} is a straightforward consequence of
\begin{eqnarray*}
P \bigl\{\psi(X_{n+1})\in B\mid\psi(X_1),\ldots,
\psi(X_n) \bigr\} &=&P \bigl\{ \psi(X_{n+1})\in B\mid
\mathcal{G}_n \bigr\}\\
&=&\alpha_n \bigl(\psi^{-1}B
\bigr)\qquad \mbox{a.s.}
\end{eqnarray*}
for each Borel set $B\subset\mathbb{R}$. This concludes the proof.
\end{pf}

When $X$ is c.i.d. (but not exchangeable) $\Vert\alpha_n-\alpha
\Vert\stackrel{\mathrm{a.s.}}\longrightarrow0$ needs not be true even if
$\alpha_c\ll
\lambda$ a.s.

\begin{ex} Let $(Z_n)$ and $(U_n)$ be independent sequences of
independent real random variables such that $Z_n\sim\mathcal
{N}(0,b_n-b_{n-1})$ and $U_n\sim\mathcal{N}(0,\break 1-b_n)$, where
$0=b_0<b_1<b_2<\cdots<1$ and $\sum_n(1-b_n)<\infty$. As shown in
Example~1.2 of \cite{BPR04},
\[
X_n=\sum_{i=1}^nZ_i+U_n
\]
is c.i.d. and $X_n\stackrel{\mathrm{a.s.}}\longrightarrow V$ for some real random
variable $V$. Since $\mu_n\stackrel{\mathrm{weak}}\longrightarrow\delta_V$ a.s.,
then $\alpha=\delta_V$ and $\alpha_c\ll\lambda$ a.s. (in fact,
$\alpha_c=0$ a.s.). However, condition \eqref{mo87xqml9} fails. In fact,
$\mathcal{L}(X_1,\ldots,X_n)\ll\lambda^n$ for all $n$, where
$\lambda$
is Lebesgue measure. Hence, $\alpha_n(\omega)\{V(\omega)\}=0$ while
$\alpha(\omega)\{V(\omega)\}=1$ for all $n$ and almost all $\omega
\in
\Omega$.
\end{ex}

We now turn to the first problem of Section \ref{intro}. Recall that
condition \eqref{nec} amounts to $\alpha_n\ll\lambda$ a.s. for all
$n\geq0$. Therefore, up to modifying $\alpha_n$ on a $P$-null set,
under condition \eqref{nec} one can write
\[
\alpha_n(\omega) (dx)=f_n(\omega,x) \lambda(dx)
\]
for each $\omega\in\Omega$, each $n\geq0$ and some function
$f_n\dvtx \Omega
\times S\rightarrow[0,\infty)$. We also~let
\begin{eqnarray*}
\mathcal{K}&=&\bigl\{K\dvtx K\mbox{ compact subset of }S\mbox{ and }\lambda (K)<\infty
\bigr\}\quad\mbox{and}
\\
 \lambda_B(\cdot)&=&\lambda(\cdot\cap B) \qquad\mbox{for all }B
\in\mathcal{B}.
\end{eqnarray*}

\begin{thmm}\label{rfg76yu8}
Suppose $X=(X_1,X_2,\ldots)$ is c.i.d. and $\mathcal{L}(X_1,\ldots
,X_n)\ll\lambda^n$ for all $n$. Then
$\alpha\ll\lambda$ a.s. if and only if, for each $K\in\mathcal{K}$,
%
\begin{eqnarray}
\begin{tabular}{p{240pt}@{}}
\label{serre123}
the sequence $(f_n(\omega,\cdot)\dvtx n\geq1
)$ is uniformly integrable,
in the space $(S,\mathcal{B},\lambda_K),$ for almost all $\omega\in\Omega.$
\end{tabular}
\end{eqnarray}
In particular, $\alpha\ll\lambda$ a.s. provided, for each $K\in
\mathcal
{K}$, there is $p>1$ such that
%
\begin{equation}
\label{mer45vk01q} \sup_n\int_K
f_n(\omega,x)^p \lambda(dx)<\infty \qquad\mbox{for almost all
}\omega\in\Omega.
\end{equation}
Moreover, for condition \eqref{mer45vk01q} to be true, it suffices that
\[
\sup_nE \biggl\{\int_K f_n^p
\,d\lambda \biggr\}<\infty.
\]
\end{thmm}

\begin{pf} If $\alpha\ll\lambda$ a.s., Theorem \ref{p09nh5t} yields
$ \Vert\alpha_n-\alpha\Vert\stackrel{\mathrm{a.s.}}\longrightarrow0$. Thus,
$f_n(\omega,\cdot)$ converges in $L_1(S,\mathcal{B},\lambda)$, for
almost all $\omega\in\Omega$, and this implies condition \eqref
{serre123}. Conversely, we now prove that $\alpha\ll\lambda$ a.s. under
condition \eqref{serre123}.

Fix a nondecreasing sequence $B_1\subset B_2\subset\cdots$ such that
$B_n\in\mathcal{B}$, $\lambda(B_n)<\infty$, and $\bigcup_nB_n=S$. Since
$\lambda(B_1)<\infty$ and $S$ is Polish, there is $K_1\in\mathcal{K}$
satisfying $K_1\subset B_1$ and $\lambda(B_1\cap K_1^c)<1$. By
induction, for each $n\geq2$, there is $K_n\in\mathcal{K}$ such that\vadjust{\goodbreak}
$K_{n-1}\subset K_n\subset B_n$ and $\lambda(B_n\cap K_n^c)<1/n$. Since
$X$ is c.i.d., condition \eqref{cideq23vg} implies
\begin{eqnarray*}
\alpha(K_m)=\lim_nE \bigl\{\alpha(K_m)\mid
\mathcal{G}_n \bigr\} =\lim_n\alpha_n(K_m)\qquad
\mbox{a.s. for all }m\geq1.
\end{eqnarray*}
Define $H=\bigcup_mK_m$ and $A_H=\{\alpha(H)=1\}$. If $\omega\in A_H$, then
\begin{eqnarray*}
\alpha(\omega) (B)=\alpha(\omega) (B\cap H)=\sup_m\alpha(\omega ) (B
\cap K_m) \qquad\mbox{for all }B\in\mathcal{B}.
\end{eqnarray*}
Moreover, $P(A_H)=1$. In fact, $\lambda(H^c)=0$ and $\alpha_n\ll
\lambda
$ a.s. for all $n$, so that
\begin{eqnarray*}
\alpha(H)=\lim_nE \bigl\{\alpha(H)\mid\mathcal{G}_n
\bigr\}=\lim_n\alpha_n(H)=1 \qquad\mbox{a.s.}
\end{eqnarray*}
Thus, to prove $\alpha\ll\lambda$ a.s., it suffices to see that
$\alpha
(\cdot\cap K_m)\ll\lambda$ a.s. for all~$m$.

Suppose \eqref{serre123} holds. Fix $m\geq1$, define $K=K_m$ and take
a set $A\in\mathcal{A}$ such that $P(A)=1$ and, for each $\omega\in A$,
\begin{eqnarray*}
&&\alpha(\omega) (K)=\lim_n\alpha_n(\omega) (K),\qquad
\alpha_n(\omega )\stackrel{\mathrm{weak}}\longrightarrow\alpha(\omega),
\\
&&\bigl(f_n(\omega,\cdot)\dvtx n\geq1\bigr)\qquad\mbox{is uniformly integrable in
}(S,\mathcal{B},\lambda_K).
\end{eqnarray*}
Let $\omega\in A$. Since $\lambda_K(S)=\lambda(K)<\infty$ and
$(f_n(\omega,\cdot)\dvtx n\geq1)$ is uniformly integrable under $\lambda_K$, there is a subsequence $(n_j)$ and a function $\psi_\omega\in
L_1(S,\mathcal{B},\lambda_K)$ such that $f_{n_j}(\omega,\cdot
)\longrightarrow\psi_\omega$ in the weak-topology of $L_1(S,\mathcal
{B},\lambda_K)$. This means that
\begin{eqnarray}
\int_{B\cap K}\psi_\omega(x) \lambda(dx)=
\lim_j\int_{B\cap
K}f_{n_j}(\omega,x)
\lambda(dx)=\lim_j\alpha_{n_j}(\omega) ({B\cap K}) \nonumber\\
\eqntext{\mbox{for
all }B\in\mathcal{B}.}
\end{eqnarray}
Therefore,
\begin{eqnarray}
\int_K\psi_\omega(x) \lambda(dx)&=&
\lim_j\alpha_{n_j}(\omega ) (K)=\alpha (\omega) (K)\quad
\mbox{and}
\nonumber\\
\int_{F\cap K}\psi_\omega(x) \lambda(dx)&=&
\lim_j\alpha_{n_j}(\omega ) (F\cap K)\leq\alpha(\omega) (F
\cap K)\nonumber \\
\eqntext{\mbox{for each closed }F\subset S.}
\end{eqnarray}
By standard arguments, the previous two relations yield
\begin{eqnarray*}
\alpha(\omega) (B\cap K)=\int_{B\cap K}\psi_\omega(x)
\lambda (dx)\qquad \mbox{for all }B\in\mathcal{B}.
\end{eqnarray*}
Thus, $\alpha(\omega)(\cdot\cap K)\ll\lambda$. This proves that
condition \eqref{serre123} implies $\alpha\ll\lambda$ a.s.

Next, since $p>1$, it is obvious that \eqref{mer45vk01q}
$\Longrightarrow$ \eqref{serre123}. Hence, it remains only to see that
condition \eqref{mer45vk01q} follows from $\sup_nE \{\int_K
f_n^p
\,d\lambda \}<\infty$.

Fix $B\in\mathcal{B}$, $p>1$, and suppose $\sup_nE \{\int_B
f_n^p
\,d\lambda \}<\infty$. Let $L_r=L_r(S,\break\mathcal{B},\lambda_B)$
for all $r$. It can be assumed $\int_B f_n(\omega,x)^p \lambda(dx)<\infty$ for
all $\omega\in\Omega$ and $n\geq1$. Thus, each $f_n\dvtx \Omega
\rightarrow
L_p$ can be seen as an $L_p$-valued random variable such that
\[
E \bigl\{\Vert f_n\Vert_{L_p} \bigr\}=E \biggl\{ \biggl(\int
_B f_n^p \,d\lambda
\biggr)^{1/p} \biggr\}\leq \biggl(E \biggl\{\int_B
f_n^p \,d\lambda \biggr\} \biggr)^{1/p}<\infty.
\]
Further, $\int f_n(\cdot,x) h(x) \lambda_B(dx)=\alpha_n(I_B h)$ is
$\mathcal{G}_n$-measurable for all $h\in L_q$, where $q=p/(p-1)$. Since
$X$ is c.i.d., condition \eqref{cideq23vg} also implies
\begin{eqnarray*}
&&E \biggl\{\int f_{n+1}(\cdot,x) h(x) \lambda_B(dx)\mid
\mathcal {G}_n \biggr\}\\
&&\qquad=E \bigl\{\alpha_{n+1}(I_B
h)\mid\mathcal{G}_n \bigr\}
\\
&&\qquad=E \bigl\{E \bigl(\alpha(I_B h)\mid\mathcal{G}_{n+1} \bigr)
\mid \mathcal{G}_n \bigr\}
\\
&&\qquad=E \bigl\{\alpha(I_B h)\mid\mathcal{G}_n \bigr\}=
\alpha_n(I_B h)
\\
&&\qquad=\int f_n(\cdot,x) h(x) \lambda_B(dx)\qquad \mbox{a.s. for all
}h\in L_q.
\end{eqnarray*}
Hence, $(f_n)$ is a $(\mathcal{G}_n)$-martingale. By Doob's maximal inequality,
\begin{eqnarray*}
E \biggl\{\sup_n\int_B f_n^p
\,d\lambda \biggr\}&=&E \Bigl\{\sup_n \Vert f_n
\Vert_{L_p}^p \Bigr\}
\\
&\leq& q^p \sup_nE \bigl\{\Vert f_n
\Vert_{L_p}^p \bigr\}=q^p \sup_nE
\biggl\{ \int_B f_n^p \,d\lambda
\biggr\}<\infty.
\end{eqnarray*}
In particular, $\sup_n\int_B f_n^p \,d\lambda<\infty$ a.s., and this
completes the proof.
\end{pf}

Some remarks on Theorem \ref{rfg76yu8} are in order.

First, for $S=[0,1]$ and a particular class of exchangeable sequences,
results similar to Theorem \ref{rfg76yu8} are in \cite{KRAFT} and
\cite{MET}.

Second,
\[
f_n(\omega,\cdot)=\frac{g_{n+1} (X_1(\omega),\ldots
,X_n(\omega),
\cdot  )}{g_n (X_1(\omega),\ldots,X_n(\omega)
)} \qquad\mbox {for almost all }\omega\in
\Omega,
\]
where each $g_n\dvtx S^n\rightarrow[0,\infty)$ is a density of $\mathcal
{L}(X_1,\ldots,X_n)$ with respect to $\lambda^n$. Thus, more
concretely, one obtains
\[
\int_K f_n^p \,d\lambda=
\frac{\int_K g_{n+1} (X_1,\ldots,X_n,
x)^p
\lambda(dx)}{g_n(X_1,\ldots,X_n)^p} \qquad\mbox{a.s.}
\]

Third, suppose $X$ exchangeable, and fix \textit{any} random probability
measure $\gamma$ on~$S$ such that $P(X\in\cdot)=\int\gamma(\omega
)^\infty(\cdot) P(d\omega)$. Then $\gamma\ll\lambda$ a.s. under the
assumptions of Theorem \ref{rfg76yu8}. In fact, $\alpha$ and $\gamma$
have the same probability distribution, when regarded as $\mathbb
{S}$-valued random variables.

A last (and important) remark deals with condition \eqref{nec}. Indeed,
even if $X$ is exchangeable, condition \eqref{nec} is not enough for
$\alpha\ll\lambda$ a.s. We close the paper showing this fact.

\begin{ex}\label{mainexamp}
Let $S=\mathbb{R}$ and $\lambda= $ Lebesgue measure. All random
variables are defined on the probability space $(\Omega,\mathcal
{A},P)$. We now exhibit an exchangeable sequence $X$ such
that $\mathcal{L}(X_1,\ldots,X_n)\ll\lambda^n$ for all $n\geq1$ and
yet $P(\alpha\ll\lambda)=0$. In fact, the support of $\alpha(\omega)$
has Hausdorff
dimension 0 for almost all $\omega\in\Omega$.

Two known facts are to be recalled. First, if $T$ and $Z$ are
independent $\mathbb{R}^n$-valued random variables, then
\[
P(T+Z\in B)=\int P(T+z\in B) P_Z(dz),
\]
where $B\in\mathcal{B}^n$ and $P_Z$ is the distribution of $Z$. Hence,
$\mathcal{L}(T+Z)\ll\lambda^n$ provided $\mathcal{L}(T)\ll\lambda^n$.
The second fact is the following:

\begin{thmm}[(Pratsiovytyi and Feshchenko)]\label{omygod47}
Let $Z_1,Z_2,\ldots$ be i.i.d. random variables with
$P(Z_1=0)=P(Z_1=1)=1/2$ and $b_1>b_2>\cdots>0$ real numbers such that
$\sum_mb_m<\infty$. Then the support of $\mathcal{L} (\sum_mb_m
Z_m )$ has Hausdorff dimension 0 whenever $\lim_m (\sum_{j>m}
b_j )^{-1} b_m=\infty$.
\end{thmm}

 Theorem \ref{omygod47} is a consequence of Theorem 8 of
\cite
{PF} (which is actually much more general).

Next, let $U_m$ and $Y_{m,n}$ be independent real random variables such that:

\begin{itemize}

\item$U_m$ is uniformly distributed on $(\frac{1}{m+1}, \frac{1}{m})$
for each $m\geq1$;

\item$P(Y_{m,n}=0)=P(Y_{m,n}=1)=\frac{1}{2}$ for all $m, n\geq1$.

\end{itemize}

 Define $V_m=U_m^m$ and
\[
X_n=\sum_{m=1}^\infty
U_m^m Y_{m,n}=\sum_{m=1}^\infty
V_m Y_{m,n}.
\]
Then, $X=(X_1,X_2,\ldots)$ is conditionally i.i.d. given $\mathcal
{V}=\sigma(V_1,V_2,\ldots)$. Precisely, for $\omega\in\Omega$ and
$B\in
\mathcal{B}$, define
\[
\alpha(\omega) (B)=P \biggl\{u\in\Omega\dvtx \sum_m
V_m(\omega) Y_{m,1}(u)\in B \biggr\}.
\]
Then, $\alpha(B)$ is a version of $P(X_1\in B\mid\mathcal{V})$ and
$P(X\in\cdot)=\int\alpha(\omega)^\infty(\cdot) P(d\omega)$. In
particular, $X$ is exchangeable. Moreover, $\mu_n\stackrel
{\mathrm{weak}}\longrightarrow\alpha$ a.s. for
\[
P (\mu_n\stackrel{\mathrm{weak}}\longrightarrow\alpha\mid\mathcal {V} )=1\qquad
\mbox{a.s.}
\]

The (topological) support of $\alpha(\omega)$ has Hausdorff dimension 0
for almost all $\omega\in\Omega$. Define in fact $b_m=V_m(\omega)$ and
$Z_m=Y_{m,1}$. By Theorem \ref{omygod47}, it suffices to verify that
%
\begin{equation}
\label{m78yuhg5} \lim_m\frac{V_m(\omega)}{\sum_{j>m} V_j(\omega)}=\infty\qquad \mbox{for almost all
}\omega\in\Omega.
\end{equation}
And condition \eqref{m78yuhg5} follows immediately from
\begin{eqnarray*}
(j+1)^{-j}&<&V_j<j^{-j} \quad\mbox{and}\\
 \sum
_{j>m} V_j&\leq&\sum_{j>m}
j^{-j}\leq\sum_{j>m} (m+1)^{-j}=
\frac{(m+1)^{-m}}{m} \qquad\mbox{a.s.}
\end{eqnarray*}

We finally prove that $\mathcal{L}(X_1,\ldots,X_n)\ll\lambda^n$ for all
$n\geq1$. Given the array $y=(y_{m,n}\dvtx m, n\geq1)$, with $y_{m,n}\in
\{
0,1\}$ for all $m, n$, define
\[
X_{n,y}=\sum_m V_m
y_{m,n}.
\]
Fix $n\geq1$ and denote $I_n$ the $n\times n$ identity matrix. If $y$ satisfies
%
\begin{eqnarray}
\label{uffa34ed5t} \pmatrix{
y_{m+1,1} &
\ldots& y_{m+1,n}
\vspace*{2pt}\cr
\ldots& \ldots& \ldots
\vspace*{2pt}\cr
y_{m+n,1} & \ldots& y_{m+n,n}}=I_n\qquad \mbox{for some }m\geq0,
\end{eqnarray}
then
\[
\begin{tabular}{p{240pt}@{}}
$(X_{1,y},\ldots,X_{n,y})=(V_{m+1},
\ldots,V_{m+n})+(R_1,\ldots,R_n)$
with $(R_1,\ldots,R_n)$ independent of $(V_{m+1},\ldots ,V_{m+n}).$
\end{tabular}
\]
In this case, since $\mathcal{L}(V_{m+1},\ldots,V_{m+n})\ll\lambda^n$,
then $\mathcal{L}(X_{1,y},\ldots,X_{n,y})\ll\lambda^n$. Hence, letting
$Y=(Y_{m,n}\dvtx m, n\geq1)$, the conditional distribution of $(X_1,\ldots,\break
X_n)$ given $Y=y$ is absolutely continuous with respect to $\lambda^n$
as far as $y$ satisfies \eqref{uffa34ed5t}. To complete the proof, it
suffices to note that
\[
P \bigl(Y=y\mbox{ for some }y\mbox{ satisfying \eqref{uffa34ed5t}} \bigr)=1.
\]

\end{ex}

%

%


\printaddresses

\end{document}